\documentclass{era-l}
\issueinfo{10}{07}{}{2004}
\dateposted{June 17, 2004}
\pagespan{58}{67}
\PII{S 1079-6762(04)00130-1}

\copyrightinfo{2004}{American Mathematical Society}
\usepackage{url}

\theoremstyle{plain}
\newtheorem{theorem}{Theorem}

\newtheorem{lemma}[theorem]{Lemma}

\newtheorem{conjecture}{Conjecture}

\numberwithin{theorem}{section}
\numberwithin{conjecture}{section}
\numberwithin{problem}{section}

\numberwithin{equation}{section}

\newcommand{\R}{{\mathbb R}}
\newcommand{\Z}{{\mathbb Z}}

\newcommand{\vol}{\mathop{\textup{vol}}}

\begin{document}

\title{The densest lattice in twenty-four dimensions}

\author{Henry Cohn}
\address{Microsoft Research,
One Microsoft Way, Redmond, WA 98052-6399}
\email{cohn@microsoft.com}

\author{Abhinav Kumar}
\address{Department of Mathematics,
Harvard University, Cambridge, MA 02138}
\email{abhinav@math.harvard.edu}
\thanks{Kumar was supported by a summer internship
in the Theory Group at Microsoft Research.}

\date{April 14, 2004}
\commby{Brian Conrey}

\begin{abstract}
In this research announcement we outline the methods used in our
recent proof that the Leech lattice is the unique densest lattice
in $\mathbb{R}^{24}$.  Complete details will appear elsewhere, but here we
illustrate our techniques by applying them to the case of lattice
packings in $\mathbb{R}^2$, and we discuss the obstacles that arise in
higher dimensions.
\end{abstract}

\maketitle

\section{Introduction}

Given any lattice $\Lambda \subset \R^n$ (i.e., a discrete
subgroup of rank $n$), one can form a sphere packing by centering
congruent balls at the points of $\Lambda$, with radius as large
as possible such that their interiors do not overlap.  Which
lattice maximizes the density of this packing, the fraction of
$\R^n$ covered by the balls?  This question arises naturally in
geometry, number theory, and information theory (sphere packings
are error-correcting codes for a continuous channel, as opposed to
a discrete channel).

The densest lattice is known for $n \le 8$.  In each of these
dimensions, it is a root lattice: $A_1$, $A_2$, $A_3$, $D_4$,
$D_5$, $E_6$, $E_7$, or $E_8$.  Furthermore each is unique up to
scaling and isometries.  The books \cite{CS} and \cite{M} are
excellent sources of information on lattice packings.  (See
\cite{Ba,Bl,G,KZ1,KZ2,V} for the original papers.)

What may be the most interesting lattices are the $E_8$ root
lattice in $\R^8$ and the Leech lattice in $\R^{24}$ (see \cite{E}
for a beautiful introduction). They are connected with many other
branches of mathematics, and are undoubtedly the densest sphere
packings in their dimensions.  However, nobody has been able to
prove their optimality.  The Leech lattice has not even been known
to be the densest lattice in its dimension, let alone the densest
sphere packing (it is possible that there are sphere packings
denser than any lattice packing).  By contrast, Blichfeldt
\cite{Bl} proved in 1935 that $E_8$ is optimal among lattices, and
in 1980 Vet\v{c}inkin \cite{V} proved its uniqueness.

In this paper we announce the following theorem:

\begin{theorem} \label{theorem:main}
The Leech lattice is the unique densest lattice in $\R^{24}$, up
to scaling and isometries of $\R^{24}$.
\end{theorem}

Our method also yields a new proof that $E_8$ is the unique
densest lattice in $\R^8$.  The proofs for $\R^8$ and $\R^{24}$
can be found in \cite{CK}.  Here, we will apply the techniques to
the case of the hexagonal lattice in $\R^2$.  That result is of no
research interest, because it is quite easy to prove using much
simpler methods than ours. Our goal is simply to illustrate our
approach in a case involving no technical difficulties.  At the
end of the paper we will outline the additional obstacles that
occur in $\R^8$ and $\R^{24}$.

Our principal technical tool will be the Poisson summation
formula.  Let $f : \R^n \to \R$ be a Schwartz function (i.e., $f$
is smooth and all its derivatives are rapidly decreasing).  We
define the Fourier transform $\widehat{f}$ of $f$ by
$$
\widehat{f}(t) = \int_{\R^n} f(x) e^{2\pi i \langle x,t \rangle}
\, dx.
$$
If $\Lambda \subset \R^n$ is a lattice, then the dual lattice
$\Lambda^*$ is defined by
$$
\Lambda^* = \{ y \in \R^n : \langle x,y \rangle \in \Z \textup{
for all $x \in \Lambda$} \}.
$$
(If $v_1,\dots,v_n$ is a basis of $\Lambda$, then the basis
$v_1^*,\dots,v_n^*$ of $\R^n$ that is dual to $v_1,\dots,v_n$
relative to the inner product forms a basis of $\Lambda^*$.)  The
Poisson summation formula states that
$$
\sum_{x \in \Lambda} f(x) = \frac{1}{\vol(\R^n/\Lambda)} \sum_{t
\in \Lambda^*} \widehat{f}(t).
$$
For a review of the standard proof, see
Appendix~\ref{appendix:Poisson}.

The most effective way to make use of a nontrivial identity such
as Poisson summation is to choose $f$ to make the two sides as
different as possible.  We will apply it using a function $f$ such
that $f(x) \le 0$ for $|x|$ sufficiently large, but
$\widehat{f}(t) \ge 0$ for all $t$.  This technique was used in
\cite{CE} to prove upper bounds on the sphere packing density (see
also \cite{C}). These bounds appear to be sharp in $\R^8$ and
$\R^{24}$, which would solve the full sphere packing problem in
those dimensions, but the apparent sharpness has not been proved.
It appears that the only dimensions in which it is sharp are $1$,
$2$, $8$, and $24$, and the one case in which it has been proved
is $\R^1$.

Our approach to proving Theorem~\ref{theorem:main} is to combine
the analytic methods from \cite{CE} with geometric and
combinatorial arguments that are specific to lattices. As
mentioned above, we will present the details here for the
two-dimensional case.  We wish to show that the obvious hexagonal
packing is the densest lattice packing of disks in $\R^2$.  For
simplicity we will choose the following normalization of the
hexagonal lattice.   Let $\Lambda_2$ be the lattice with basis
$v,w$, where
$$
v = \begin{pmatrix} \sqrt{2}\\ 0 \end{pmatrix}
$$
and
$$
w = \begin{pmatrix} \sqrt{2}/2\\
\sqrt{6}/2\end{pmatrix}.
$$
The advantage of this normalization is that $\Lambda_2$ is an even
integral lattice.  In other words, all inner products of lattice
points are integers, and furthermore every lattice point has even
norm.  Here, $|v|^2=|w|^2=2$ and $\langle v,w \rangle = 1$.

Let $\Lambda \subset \R^2$ be any lattice that is at least as
dense as $\Lambda_2$.  We scale $\Lambda$ so that
$\vol(\R^2/\Lambda)=1$ (we refer to $\vol(\R^2/\Lambda)$ as the
covolume of $\Lambda$). Unfortunately, $\vol(\R^2/\Lambda_2) =
\sqrt{3}$, so the scalings of $\Lambda$ and $\Lambda_2$ are not
compatible. That will cause us no problems, and in fact each
scaling will prove convenient at a certain point in the proof.

Our overall plan is as follows:
\begin{enumerate}
\item Use Poisson summation to study the short
vectors in $\Lambda$.

\item Apply this information to see that $\Lambda$ is, up to scaling,
a small perturbation of $\Lambda_2$.

\item Prove that $\Lambda_2$ is a strict local optimum for the packing density.

\item Conclude that either $\Lambda$ is the same as $\Lambda_2$ (up to
scaling and isometries) or it is strictly less dense, which would
contradict our assumptions.
\end{enumerate}

Before beginning, it is worth noting that the density of $\Lambda$
simply depends on its minimal vector length (i.e., the length of
the shortest nonzero vectors in $\Lambda$).  Because
$\vol(\R^2/\Lambda) = 1$, there is one sphere per unit volume in
space, and the density of the packing equals the volume of such a
sphere.  More generally, the density of any lattice $L \subset
\R^n$ equals
$$
\frac{\vol(B_r)}{\vol(\R^n/L)},
$$
where $r$ denotes the packing radius and $B_r$ denotes a ball of
radius $r$. The packing radius is half of the minimal vector
length, because that is the radius at which adjacent spheres are
tangent.  Thus, our goal becomes understanding the minimal vector
length of $\Lambda$.

Because $\Lambda$ is at least as dense as $\Lambda_2$, its minimal
vectors have length at least $(4/3)^{1/4}$.  To see why, recall
that $\vol(\R^2/\Lambda_2) = \sqrt{3}$ and its minimal vectors
have length $\sqrt{2}$. Rescaling $\Lambda_2$ by a factor of
$3^{-1/4}$ makes the covolume $1$ and yields minimal vector length
$(4/3)^{1/4}$.

\section{Applications of Poisson summation}

Define $f : \R^2 \to \R$ by $f(x) = p_f(2\pi|x|^2) e^{-\pi|x|^2}$,
where
$$
p_f(u) =  20812+756u+1107u^2-216u^3.
$$
Then one can calculate that $\widehat{f}(t) =
p_{\widehat{f}}(2\pi|t|^2) e^{-\pi|t|^2}$, where
$$
p_{\widehat{f}}(u) = 20812+5940u-2781u^2+216u^3 =
(43+24u)(-22+3u)^2.
$$
We have $f(0) = \widehat{f}(0)$, $f(x) < 0$ for $|x| \ge 1.084$,
and $\widehat{f}(t) \ge 0$ for all $t$.  (Terminating decimal
expansions such as $1.084$ represent exact rational numbers, not
floating point approximations.)  Up to scaling, $p_f$ is the
unique polynomial of degree $3$ such that $f(0)=\widehat{f}(0)$
and $p_{\widehat{f}}$ has a double root at $22/3$, and that is how
it was constructed.

\begin{lemma}
$\Lambda$ contains a nonzero vector of length at most $1.084$.
\end{lemma}

Note that our lower bound $(4/3)^{1/4}$ for the minimal vector
length is between $1.074$ and $1.075$.

\begin{proof}
By Poisson summation,
$$
\sum_{x \in \Lambda} f(x) = \sum_{t \in \Lambda^*} \widehat{f}(t)
\ge \widehat{f}(0),
$$
because $\widehat{f}(t) \ge 0$ for all $t$. On the other hand,
$f(0)=\widehat{f}(0)$ and $f(x) < 0$ for $|x| \ge 1.084$, so it
follows that
$$
\sum_{x \in \Lambda,\ 0 < |x| < 1.084} f(x) > 0.
$$
Therefore the minimal vector length of $\Lambda$ is at most
$1.084$.
\end{proof}

Call a vector in $\Lambda$ ``nearly minimal'' if it has length in
$[(4/3)^{1/4},1.114)$ (the reason for the upper bound of 1.114
will become clear below).

\begin{lemma} \label{lemma:atmostsix}
There are at most six nearly minimal vectors in $\Lambda$.
\end{lemma}

\begin{proof}
If $x$ and $y$ are nearly minimal and $\theta \in [0,\pi]$ is the
angle between them, then
$$
\cos \theta = \frac{|x|^2+|y|^2-|x-y|^2}{2|x||y|} \le \frac{2\cdot
1.114^2-(4/3)^{1/2}}{2\cdot (4/3)^{1/2}} < 0.575.
$$
Thus, $\theta > \cos^{-1}0.575 > 2\pi \cdot 0.152$.  There is no
room for seven vectors to be separated by such an angle: if we
assign to each nearly minimal vector the arc on the unit circle
consisting of all points within angle $2\pi \cdot 0.076$ of it,
then these arcs do not overlap.  Because each arc has length $2\pi
\cdot 0.152$, if there were seven of them then the arc length of
the unit circle would be at least $2\pi \cdot 1.064$.
\end{proof}

In the other direction, we can use Poisson summation to see that
there are more than five nearly minimal vectors in $\Lambda$ (so
there must be exactly six).  First we need one more lemma:

\begin{lemma}
Every nonzero vector in $\Lambda$ is either nearly minimal or else
has length at least $1.62$.
\end{lemma}

\begin{proof}
First, note that $f(x)$ is a decreasing function of $|x|$
for $|x| \in [0,1.084]$.  In the inequality
$$
\sum_{x \in \Lambda} f(x) \ge \widehat{f}(0),
$$
Lemma~\ref{lemma:atmostsix} implies that there are at most six
positive terms on the left hand side other than $x=0$.  Each of
them is at most $f\big((4/3)^{1/4}\big)$, where for $r \in
[0,\infty)$ we write $f(r)$ to indicate the common value $f(x)$
with $|x|=r$. Thus,
\begin{equation}\label{ineq}
6 f\big((4/3)^{1/4}\big) + \sum_{x \in \Lambda,\ f(x)<0} f(x) \ge
0.
\end{equation}
One can check that $f(x) < -3f\big((4/3)^{1/4}\big) < 0$ for $|x|
\in [1.114, 1.62]$, and hence there is no $x \in \Lambda$
satisfying $|x| \in [1.114, 1.62]$ (because the $x$ and $-x$ terms
in \eqref{ineq} would combine to make the left hand side
negative).
\end{proof}

\begin{lemma}
There are more than five nearly minimal vectors in $\Lambda$.
\end{lemma}

\begin{proof}
Define $g : \R^2 \to \R$ by $g(x) = p_g(2\pi|x|^2) e^{-\pi|x|^2}$,
where
$$
p_g(u) =  (13-u)(1075+220u+69u^2).
$$
Then $\widehat{g}(t) = p_{\widehat{g}}(2\pi|t|^2) e^{-\pi|t|^2}$,
where
$$
p_{\widehat{g}}(u) = (401+69u)(u-7)^2.
$$
We have $g(x) \le 0$ for $|x| \ge 1.62$, and $\widehat{g}(t) \ge
0$ for all $t$.

It follows from Poisson summation that if $\mathcal{M}$ denotes
the set of nearly minimal vectors in $\Lambda$, then
$$
g(0) + \sum_{x \in \mathcal{M}} g(x) \ge \widehat{g}(0).
$$
On the other hand, $g(x)$ is a decreasing function of $|x|$ for
$|x| \in [(4/3)^{1/4},1.114)$, so
$$
g(0) + |\mathcal{M}| g\big((4/3)^{1/4}\big) \ge \widehat{g}(0).
$$
Thus,
$$
|\mathcal{M}| \ge
\frac{\widehat{g}(0)-g(0)}{g\big((4/3)^{1/4}\big)} > 5.89,
$$
as desired.
\end{proof}

\section{Arrangement of the nearly minimal vectors}

There are two distinct nearly minimal vectors $x,y \in \Lambda$
within angle $2\pi/6$ of each other, because that is the average
angle between them as one moves around the circle. Then $x-y$ must
be nearly minimal as well: if $\theta$ is the angle between them,
then
$$
|x-y|^2 = |x|^2+|y|^2-2|x||y|\cos \theta \le |x|^2+|y|^2-|x||y|
$$
because $\theta \le 2\pi/6$, so
$$
|x-y|^2 \le 2\cdot 1.114^2 - (4/3)^{1/2} < 1.33.
$$
It follows that $x$, $y$, $x-y$, and their negatives form the full
list of nearly minimal vectors.

\begin{lemma}
The vectors $x$ and $y$ form a basis of $\Lambda$.
\end{lemma}

\begin{proof}
The nearly minimal vectors are all in the span of $x$ and $y$, so
we need only prove that the nearly minimal vectors span $\Lambda$.
Suppose not, and that $z \in \Lambda$ is the smallest vector not
in their span. There cannot exist a nearly minimal vector $u$ such
that $|u-z| < |z|$.  For every nearly minimal $u$, the angle
$\theta$ between $u$ and $z$ satisfies
$$
\cos \theta = \frac{|u|^2+|z|^2-|u-z|^2}{2|u||z|} \le
\frac{|u|}{2|z|} < \frac{1}{2} = \cos \frac{2\pi}{6}
$$
because $|u-z| \ge |z|$ and $|u|<|z|$.  In other words, no nearly
minimal vector is within angle $2\pi/6$ of $z$.  That is
impossible, because then $z$ together with the nearly minimal
vectors would form seven vectors separated by angles of at least
$2\pi \cdot 0.152$.
\end{proof}

Thus, we have shown that $\Lambda$ has a basis $x,y$ such that
$|x|,|y| \in [(4/3)^{1/4},1.114)$ and the angle between $x$ and
$y$ is in $[2\pi \cdot 0.152, 2\pi/6]$.

Consider the lattice $3^{1/4}\Lambda$, which we have rescaled so
that its covolume is $3^{1/2}$, the same as that of $\Lambda_2$.
The rescaled basis vectors $3^{1/4}x,3^{1/4}y$ have lengths in
$[\sqrt{2},1.467)$ and therefore norms in $[2,2.16)$. If $\theta$
denotes the angle between them then
$$
\langle 3^{1/4}x,3^{1/4}y \rangle = 3^{1/2}|x||y|\cos\theta \in
[2\cos (2\pi/6), 1.467^2 \cos (2\pi \cdot 0.152)] \subset
[1,1.243].
$$
Therefore the Gram matrix of $3^{1/4}\Lambda$ with respect to the
basis $3^{1/4}x,3^{1/4}y$ has entries within $0.243$ of those of
the Gram matrix of $\Lambda_2$ with respect to $v,w$.  It is in
this sense that $\Lambda$ is a perturbation of $\Lambda_2$.

\section{Local optimality}

In this section, we prove that $\Lambda_2$ is a strict local
optimum for density.  It is more convenient to deal with quadratic
forms than with lattices.  Let $Q$ be the quadratic form
corresponding to $\Lambda_2$ and its basis $v,w$ (i.e., $Q(s,t) =
|sv+tw|^2$ for $s,t \in \R$). The Gram matrix of $\Lambda_2$ with
respect to the basis $v,w$ is
$$
\begin{pmatrix}
2 & 1\\
1 & 2
\end{pmatrix},
$$
and we have
$$
Q(s,t) = \begin{pmatrix}s & t\end{pmatrix}
\begin{pmatrix}
2 & 1\\
1 & 2\end{pmatrix} \begin{pmatrix}s\\ t\end{pmatrix}.
$$
Note that the Gram matrix of a lattice basis determines the
lattice up to isometries, so we lose no information by focusing on
it.

The minimal norm $M$ of $Q$ is the minimum of $Q(s,t)$ over $(s,t)
\in \Z^2\setminus \{(0,0)\}$ (i.e., $M=2$), and the determinant
$D$ is defined by
$$
D = \det \begin{pmatrix}
2 & 1\\
1 & 2
\end{pmatrix} = 3.
$$
In these terms, the density of $\Lambda_2$ equals
$(\pi/4)D^{-1/2}M$, because $D^{1/2}$ is the covolume
$\vol(\R^2/\Lambda_2)$ and $M$ is the square of the minimal vector
length of $\Lambda_2$.

We will show that if $Q$ is slightly perturbed, then either the
perturbation is itself proportional to $Q$ (which corresponds to
rescaling $\Lambda_2$) or else the density strictly decreases. We
use symmetric perturbations because Gram matrices are symmetric.

Let $\max \{|a|,|b|,|c|\} = \rho>0$, and let $Q_\rho$ be the
perturbation of $Q$ with matrix
$$
\begin{pmatrix}
2+a & 1+b\\
1+b & 2+c
\end{pmatrix}.
$$
Let
$$
D_\rho = \det \begin{pmatrix}
2+a & 1+b\\
1+b & 2+c
\end{pmatrix} = 3 + 2(a+c-b) + (ac-b^2),
$$
and let $M_\rho$ be the minimal norm of $Q_\rho$.

\begin{lemma}
If $\rho < 12/47$, then either $D_\rho^{-1/2}M_\rho < D^{-1/2}M$
or $Q_\rho$ is proportional to $Q$.
\end{lemma}

\begin{proof}
We begin by assuming that $a+c=b$.  In other words, the linear
terms in the expansion of $D_\rho$ cancel.  The minimal norm of
$Q$ occurs at $(1,0)$, $(0,1)$, $(1,-1)$, and their negatives, and
at these points $Q_\rho$ takes on the values $2+a$, $2+c$, and
$2+a-2b+c$, respectively.  Consider the perturbations $a$, $c$,
and $a-2b+c$ of these values away from $2$.  It follows from
$a+c=b$ that $a-2b+c=-b$ and that these three perturbations sum to
$0$.  Thus, at least one of them must be negative (if all three
vanished then $\rho$ would too). In fact, one of them is at most
$-\rho/2$: if $a$, $c$, or $-b$ equals $-\rho$ then that is
trivial, and otherwise one of them equals $\rho$ and then one of
the other two is at most $-\rho/2$ because they sum to $0$.

Thus, if $a+c=b$, then $M_\rho \le 2 - \rho/2$.  When we combine
that with
$$
D_\rho = 3 + (ac-b^2) \ge 3 - 2\rho^2,
$$
we find that
$$
D_\rho^{-1/2}M_\rho \le (3-2\rho^2)^{-1/2}(2-\rho/2).
$$
When $0 < \rho < 24/35$,
$$
(3-2\rho^2)^{-1/2}(2-\rho/2) < 3^{-1/2}\cdot 2 = D^{-1/2}M,
$$
so the density of the perturbed lattice is strictly less than that
of $\Lambda_2$.

In general we cannot assume $a+c=b$.  However, $Q_\rho$ is
proportional to a perturbation of $Q$ in which that equation
holds.  Set $A = (a-2c+2b)/(3+a+c-b)$, $B = (4b-a-c)/(3+a+c-b)$,
and $C = (2b+c-2a)/(3+a+c-b)$.  Then
$$
\begin{pmatrix}
2+a & 1+b\\
1+b & 2+c
\end{pmatrix}=
\left(1+(a+c-b)/3\right) \begin{pmatrix}
2+A & 1+B\\
1+B & 2+C
\end{pmatrix}
$$
and $A+C=B$. Clearly, $|A| \le 5\rho/(3-3\rho)$, $|B| \le
2\rho/(1-\rho)$, and $|C| \le 5\rho/(3-3\rho)$.  If $\rho <
12/47$, then $|A|,|B|,|C| < 24/35$, and so the perturbed lattice
is either itself hexagonal or else strictly less dense than
$\Lambda_2$.
\end{proof}

\begin{theorem} \label{theorem:toy}
The hexagonal lattice $\Lambda_2$ is the unique densest lattice in
$\R^2$, up to scaling and isometries of $\R^2$.
\end{theorem}

\begin{proof}
We have seen that if the Gram matrix of $\Lambda_2$ is perturbed
by at most $12/47$, then either the density strictly decreases or
the lattice remains the same (up to a similarity).  On the other
hand, the Gram matrix of $3^{1/4}\Lambda$ with respect to
$3^{1/4}x,3^{1/4}y$ has entries within $0.243$ of those of the
Gram matrix of $\Lambda_2$.  Observing that $0.243 < 12/47$
completes the proof.
\end{proof}

\section{Relationship with the higher-dimensional cases}

The techniques we apply in \cite{CK} are completely analogous to
those we have used here, and each of the lemmas from this paper
has a counterpart there.  However, there are a couple of key steps
in \cite{CK} that are not apparent from the two-dimensional case,
as well as some technical obstacles that must be overcome.

\subsection{Spherical codes}

The biggest conceptual difference between this proof and the one
in \cite{CK} is in how the arrangement of nearly minimal vectors
is studied.  In each case, when rescaled to the unit sphere they
form a spherical code, i.e., a subset of the sphere that contains
no pair of points closer than some given angle.  In $\R^2$ the
unit sphere is a circle, and distributing points on a circle with
at least a given angle between them is trivial.  In higher
dimensions this problem is much more subtle.

Independently, Levenshtein \cite{L} and Odlyzko and Sloane
\cite{OS} used linear programming bounds for spherical codes to
solve the kissing problem in $\R^8$ and $\R^{24}$: how many unit
balls can be placed tangent to a given one, if they may not
overlap except tangentially?  The answers are $240$ and $196560$,
respectively, and the arrangements come from the $E_8$ and Leech
lattice packings.

In \cite{CK} we use linear programming bounds in a similar way to
bound the number of nearly minimal vectors and the angles between
them.

\subsection{Spherical designs and association schemes}

It is not enough in higher dimensions simply to compute the number
of nearly minimal vectors.  To determine the configuration more
precisely, we apply techniques derived from \cite{DGS}.  In
particular, we formulate a notion of an approximate spherical
design, and prove that the nearly minimal vectors form one.  We
then show that if they are grouped according to the (approximate)
angles between them, then they form an association scheme.
Finally, we prove that the association scheme with these
parameters is unique, which lets us conclude that the nearly
minimal vectors form a perturbation of the desired configuration.
In principle all these results have analogues in $\R^2$, with
similar proofs, but we did not require them for
Theorem~\ref{theorem:toy}, and including them would have
substantially lengthened this article.

\subsection{Local optimality and angle bounds}

The proof of strict local optimality is analogous to that given
here. We use Voronoi's theorem characterizing locally optimal
lattices (as proved in \cite[\S39]{GL}).  However, it is more
difficult to supply numerical bounds on how large the
perturbations can be without letting the density increase.
Section~10 of \cite{CK} is devoted to computing such a bound.

In the other direction, it is also more difficult to prove bounds
on how close the unknown lattice is to $E_8$ or the Leech lattice.
Section~7 applies the uniqueness of the association scheme to this
problem, in a way that has no analogue in $\R^2$ because fewer
possible angles occur between minimal vectors in $\Lambda_2$ than
in the Leech lattice.

\subsection{Computer calculations}

Completing the proofs in $\R^8$ and $\R^{24}$ requires much
sharper estimates than in $\R^2$.  For example, along the way we
prove that no sphere packing in $\R^{24}$ can exceed the Leech
lattice's density by a factor of $1 + 1.65 \cdot 10^{-30}$.
Proving such bounds using these techniques requires applying
Poisson summation to much more complicated functions.  For
example, the analogue in \cite{CK} of the polynomial $p_f$ from
this paper has degree $803$ and rational coefficients with
denominator $10^{3000}$.  Checking the desired properties clearly
requires a computer.  We have arranged the calculations so as to
run in one hour on an ordinary personal computer, and we have made
our code available (see Appendix~A in \cite{CK}).  All of our
calculations use exact rational arithmetic and are completely
rigorous.

In addition to the difficulty of checking the assertions about
$f$, there is the issue of how to construct this function. We used
a computer to search for polynomials with the desired properties,
and optimized them using a high-dimensional version of Newton's
method.  The whole process took our computers approximately a
month, but fortunately one can check the proof simply by verifying
the final answer, with no need to reconstruct it from scratch.

\section{Future prospects}

There is no likelihood that our techniques will work in any
dimension except $1$, $2$, $8$, or $24$.  They depend on using the
methods of \cite{CE} to prove nearly sharp upper bounds for the
packing density, and that does not seem to happen in any other
dimension. It seems counterintuitive that determining the densest
lattice in $\R^{16}$ appears much more difficult than in $\R^8$ or
$\R^{24}$, but that simply seems to be the case.

One might hope to solve the full sphere packing problem using our
techniques.  We cannot absolutely rule that out, but we do not
consider it feasible.  There are at least two substantial
obstacles:
\begin{enumerate}
\item It is not known whether $E_8$ or the Leech lattice is even locally
optimal when the perturbations leave the space of lattices.

\item In a general sphere packing, different spheres can be tangent
to different numbers of other spheres (unlike the case of
lattices).  Our method of proving lower bounds for the number of
nearly minimal vectors can be used to bound the \textit{average\/}
number of near neighbors of a sphere in a dense packing.  However,
we cannot rule out the possibility of a small minority of spheres
with few neighbors.
\end{enumerate}

We do not know how to deal with the first problem, but perhaps it
could be done.  The second problem appears quite fundamental, and
we suspect that dealing with it would require a major advance in
the theory of sphere packing.

Instead, we are convinced that the right way to solve the sphere
packing problem in $\R^8$ and $\R^{24}$ is to prove the following
conjecture from \cite{CE}.  Call a continuous function $f: \R^n
\to \R$ admissible if there is a constant $\delta>0$ such that
both $|f(x)|$ and $|\widehat{f}(x)|$ are bounded above by a
constant times $(1+|x|)^{-n-\delta}$.  (Note that the bound on $f$
implies that the integral defining $\widehat{f}$ converges.)

\begin{conjecture} \label{conjecture:cohnelkies}
For $(n,r) = (8,\sqrt{2})$ or $(24,2)$, there is an admissible
function $f : \R^n \to \R$ such that $f(0) = \widehat{f}(0) = 1$,
$f(x) \le 0$ for $|x| \ge r$, and $\widehat{f}(t) \ge 0$ for all
$t$.
\end{conjecture}

If Conjecture~\ref{conjecture:cohnelkies} holds, then Theorem~3.2
of \cite{CE} implies that $E_8$ and the Leech lattice are the
densest sphere packings in $\R^8$ and $\R^{24}$, respectively. In
\cite{CK} we achieve $r \le 2\big(1+6.851 \cdot 10^{-32}\big)$
when $n=24$, which is numerical evidence that $r=2$ could be
achieved.

\section*{Acknowledgements}

We thank Amanda Beeson and James Bernhard for comments on our paper.

\appendix
\section{Poisson summation}
\label{appendix:Poisson}

To make this paper self-contained, we provide here the standard
proof of the Poisson summation formula:

\begin{theorem}[Poisson summation]
Let $f : \R^n \to \R$ be a Schwartz function and $\Lambda \subset
\R^n$ a lattice.  Then
$$
\sum_{x \in \Lambda} f(x) = \frac{1}{\vol(\R^n/\Lambda)} \sum_{t
\in \Lambda^*} \widehat{f}(t).
$$
\end{theorem}

\begin{proof}
It is easier to prove a more general formula. Define $F : \R^n \to
\R$ by
$$
F(z) = \sum_{x \in \Lambda} f(x+z).
$$
In other words, $F$ is $f$ made periodic modulo $\Lambda$, and we
wish to compute $F(0)$. Because $f$ is a Schwartz function, the
sum defining $F$ converges and defines a $C^\infty$ function.  It
follows that we can expand $F$ as a Fourier series, in particular
as a linear combination of the exponential functions that are
periodic modulo $\Lambda$.  The function
$$
z \mapsto e^{-2\pi i \langle t,z \rangle}
$$
is periodic modulo $\Lambda$ iff $t \in \Lambda^*$.  Thus,
$$
F(z) = \sum_{t \in \Lambda^*} c_t e^{-2\pi i \langle t,z \rangle}
$$
for some coefficients $c_t$.  We can compute $c_t$ using
orthogonality: if $D$ is a fundamental domain for $\Lambda$, then
$$
c_t = \frac{1}{\vol(D)} \int_D F(z) e^{2\pi i \langle t,z \rangle}
\, dz.
$$
When we substitute $\vol(D) = \vol(\R^n/\Lambda)$ and the
definition of $F(z)$ into this equation, we find that
\begin{eqnarray*}
c_t &=& \frac{1}{\vol(\R^n/\Lambda)} \sum_{x \in \Lambda} \int_D
f(x+z) e^{2\pi i \langle t,z \rangle} \, dz\\
&=& \frac{1}{\vol(\R^n/\Lambda)} \sum_{x \in \Lambda} \int_{D+x}
f(z) e^{2\pi i \langle t,z \rangle} \, dz\\
&=& \frac{1}{\vol(\R^n/\Lambda)} \int_{\R^n} f(z) e^{2\pi i
\langle t,z \rangle} \, dz\\
&=& \frac{1}{\vol(\R^n/\Lambda)} \widehat{f}(t).
\end{eqnarray*}
Here $D+x = \{d+x : d \in D\}$, and these sets tile $\R^n$ as $x$
ranges over $\Lambda$.

It follows that
$$
\sum_{x \in \Lambda} f(x+z) = \frac{1}{\vol(\R^n/\Lambda)}\sum_{t
\in \Lambda^*} \widehat{f}(t) e^{-2\pi i \langle t,z \rangle}.
$$
Setting $z=0$ yields the desired result.
\end{proof}

\end{document}